\documentclass[12pt,leqno]{article}

\usepackage{amsthm}

\usepackage{amsmath}
\usepackage{amscd}
\usepackage{amsmath}
\usepackage{amssymb}
\usepackage{latexsym}
\makeatletter

\@addtoreset{equation}{section}
\makeatother

\newtheorem{thm}{Theorem}[section]
\newtheorem{pr}[thm]{Proposition}
\newtheorem{df}[thm]{Definition}
\newtheorem{lm}[thm]{Lemma}
\newtheorem{cor}[thm]{Corollary}

\input xy \xyoption{all} \CompileMatrices

\begin{document}

\title{Characteristic cycle
of the exterior product\\
of constructible sheaves}
\author{Takeshi Saito}

\maketitle

\begin{abstract}
We show that
the characteristic cycle 
of the exterior product
of constructible complexes
is the exterior product
of the characteristic cycles
of factors.
This implies the compatibility
of characteristic cycles
with smooth pull-back
which is a first step in the proof
of the index formula.
\end{abstract}

The characteristic cycle of
a constructible complex on
a smooth scheme over
a perfect field 
is defined as a cycle on the cotangent
bundle \cite{CC}
supported on the singular support
\cite{Be}.
It is characterized by
the Milnor formula \cite[(5.15)]{CC}
for the vanishing cycles
defined for morphisms to curves.

We prove
a formula (\ref{eqCCpr}) for the external product
in Theorem \ref{thmpr}.
Theorem \ref{thmpr}
implies the compatibility of
characteristic cycles with
smooth pull-back Corollary \ref{corpull}
which is a first step in the proof of
the index formula \cite[Theorem 7.13]{CC}.
Note that Theorem \ref{thmpr}
is proved without using
the results in \cite{CC}
after Proposition 5.17 loc.\! cit.\!
included.
Corollary \ref{corpr} corresponds
to \cite[Corollary 5.4.14]{KSc}.

We briefly sketch the idea of proof
of Theorem \ref{thmpr}.
First we show that
the external product
is micro-supported
on the external product
of the singular supports
of the factors.
We deduce this from projection
formulas for nearby cycles
over general base schemes
in \cite{W} recalled in Section \ref{spsi}.
The formula (\ref{eqCCpr})
for characteristic cycle
is deduced from
the Thom-Sebastiani
formula \cite{TS} and
a conductor formula (\ref{eqa*}) for
the additive convolution
\cite[Corollary 5.12]{TS}
which is an analogue for torsion
coefficient of \cite[Proposition (2.7.2.1)]{La}.

Corollary \ref{corpull}
of Theorem \ref{thmpr}
is a first step of the proof of
index formula in general dimension 
\cite[Theorem 7.13]{CC}.
The formula (\ref{eqa*}) for convolution
on which the proof of
Theorem \ref{thmpr} is based
is essentially equivalent to
the formula (\ref{eqCCpr})
in the case where $\dim X=\dim Y=1$.
This can be deduced from
the special case of the index formula
in dimension 2 proved earlier in 
\cite[Theorem 3.19]{surface}.

The author thanks Luc Illusie
for discussion on nearby cycles
over general base schemes,
the additive convolution etc.
He thanks to Zheng Weizhe
for an unpublished notes
\cite{W}.
The research was partially supported
by JSPS Grants-in-Aid 
for Scientific Research
(A) 26247002.


\section{Nearby cycles
and projection formulas}\label{spsi}

Let $f\colon X\to S$ be
a morphism of schemes.
For the definition and properties of
the vanishing topos
$X\overset\gets\times_SS$
and the diagram
\begin{equation}
\xymatrix{
&X\ar[dl]_{\rm id}
\ar[d]^{\Psi_f}\ar[dr]^f\\
X\ar[dr]_f
&X\overset\gets\times_SS
\ar[l]_{p_1}\ar[r]^{p_2}&
S\ar[dl]^{\rm id}
\\
&S,
}\label{eqXXS}
\end{equation}
we refer to
\cite{Il}, \cite{TS}, \cite{Or}.

Assume that $f\colon X\to S$ is
a morphism of finite type of
noetherian schemes.
Let $\Lambda$ be
a finite field of characteristic $\ell$
invertible on $S$
and let ${\cal F}$
and ${\cal G}$
be complexes bounded above of 
$\Lambda$-modules
on $X$ and on $S$
respectively.
A canonical morphism
\begin{equation}
R\Psi_f{\cal F}\otimes^L_\Lambda
p_2^*{\cal G}
\to 
R\Psi_f({\cal F}\otimes^L_\Lambda
f^*{\cal G})
\label{eqPs}
\end{equation}
on $X\overset\gets\times_SS$
is defined as the adjoint
of
$\Psi_f^*(R\Psi_f{\cal F}\otimes^L_\Lambda
p_2^*{\cal G})
=
\Psi_f^*R\Psi_f{\cal F}\otimes^L_\Lambda
f^*{\cal G}
\to 
{\cal F}\otimes^L_\Lambda
f^*{\cal G}$
induced by the adjunction
$\Psi_f^*R\Psi_f{\cal F}
\to {\cal F}$,
since $\Psi_f$ is of finite cohomological
dimension by \cite[Proposition 3.1]{Or}.

\begin{lm}[{\rm \cite[Proposition 4]{W}}]\label{lmpf1}
Let $f\colon X\to S$ be
a morphism of 
finite type of noetherian schemes
and let ${\cal F}$
and ${\cal G}$
be complexes bounded above of 
$\Lambda$-modules
on $X$ and on $S$
respectively.
We assume that
the formation of
$R\Psi_f{\cal F}$
commutes with finite base change.
Then,
the canonical morphism
\begin{equation*}
R\Psi_f{\cal F}\otimes^L_\Lambda
p_2^*{\cal G}
\to 
R\Psi_f({\cal F}\otimes^L_\Lambda
f^*{\cal G})
\leqno{\rm (\ref{eqPs})}
\end{equation*}
on $X\overset\gets\times_SS$
is an isomorphism.
\end{lm}

Further, let
$h\colon W\to X$
be a morphism of schemes
and consider the commutative
diagram
\begin{equation}
\xymatrix{
&
X
\\
W\overset\gets\times_XX
\ar[r]^{\overset\gets h}
\ar[d]_{\overset\gets f}
\ar[ru]^{p_2}
&
X\overset\gets\times_XX
\ar[u]^{p_2}
\ar[d]^{\overset\gets f}
&
X\ar[l]_{\Psi_{\rm id}}\ar[ld]^{\Psi_f}
\ar[lu]_{\rm id}
\\
W\overset\gets\times_SS
\ar[r]^{\overset\gets h}&
X\overset\gets\times_SS}
\label{eqwxs}
\end{equation}
of vanishing toposes.
By \cite[Proposition 3.1]{Or},
$\overset\gets f_*$
is of finite cohomological dimension.
Let ${\cal F}$
and ${\cal G}$ be
complexes of $\Lambda$-modules
on $X$ and on $W$
respectively.
We define a base change morphism
\begin{equation}
\overset\gets h^*
R\Psi_f
{\cal F}
\otimes^L_\Lambda
p_1^*{\cal G}
\to
R\overset\gets f_*(
p_2^*{\cal F}
\otimes^L_\Lambda
p_1^*{\cal G})
\label{eqcan}
\end{equation}
on $W\overset\gets\times_SS$
as the adjoint of
the morphism
$\overset\gets f^*(
\overset\gets h^*
R\Psi_f
{\cal F}
\otimes^L_\Lambda
p_1^*{\cal G})
=
\overset\gets f^*
\overset\gets h^*
R\Psi_f
{\cal F}
\otimes^L_\Lambda
p_1^*{\cal G}
\to
p_2^*{\cal F}
\otimes^L_\Lambda
p_1^*{\cal G}$
defined as follows.
We identify
$
\overset\gets f^*
\overset\gets h^*
R\Psi_f
=
\overset\gets h^*
\overset\gets f^*
R\overset\gets f_*
R\Psi_{\rm id}$
and
$p_2^*
=
\overset\gets h^*
p_2^*=
\overset\gets h^*
R\Psi_{\rm id}$
by the isomorphism
$p_2^*\to
R\Psi_{\rm id}$
\cite[Proposition 4.7]{Il}
defined as the adjoint of
$\Psi_{\rm id}^*p_2^*\to
{\rm id}$.
Then, the morphism in question
is induced by the adjunction
$\overset\gets f^*
R\overset\gets f_*
\to {\rm id}.$

\begin{lm}[{\rm \cite[Proposition 5]{W}}]\label{lmpf2}
Let $f\colon X\to S$ be
a morphism of 
finite type of noetherian schemes
and $h\colon W\to X$
be a morphism of schemes.
Let ${\cal F}$
and ${\cal G}$
be complexes bounded above of 
$\Lambda$-modules
on $X$ and on $W$
respectively.
We assume that
the formation of
$R\Psi_f{\cal F}$
commutes with finite base change.
Then,
the canonical morphism
\begin{equation*}
\overset\gets h^*
R\Psi_f
{\cal F}
\otimes^L_\Lambda
p_1^*{\cal G}
\to
R\overset\gets f_*(
p_2^*{\cal F}
\otimes^L_\Lambda
p_1^*{\cal G})
\leqno{\rm (\ref{eqcan})}
\end{equation*}
on $W\overset\gets\times_SS$
is an isomorphism.
\end{lm}

We recall an interpretation
of local acyclicity
in terms of vanishing topos.

\begin{pr}[{\rm \cite[Proposition 1.7]{CC}}]\label{lmapp}
Let $f\colon X\to S$
be a morphism of schemes.
Then, for a complex ${\cal F}\in D^+(X)$
bounded below, 
the following conditions {\rm (1)} and {\rm (2)}
are equivalent:

{\rm (1)}
The morphism $f\colon X\to S$ is 
locally acyclic relatively to ${\cal F}$.

{\rm (2)}
The formation of
$R\Psi_f{\cal F}$
commutes with every
finite base change $T\to S$
and 
the canonical morphism
$p_1^*{\cal F}\to R\Psi_f{\cal F}$
is an isomorphism.

{\rm (3)}
The canonical morphism
$p_1^*{\cal F}_T\to
R\Psi_{f_T}{\cal F}_T$
is an isomorphism
for every finite morphism $T\to S$,
the cartesian diagram
$$\begin{CD}
X@<<< X_T\\
@VfVV @VV{f_T}V\\
S@<<< T
\end{CD}$$
and the pull-back ${\cal F}_T$
of ${\cal F}$ on $X_T$.
\end{pr}

\begin{cor}\label{corpf}
Let $f\colon X\to S$ be
a morphism of 
finite type of noetherian schemes
and let ${\cal F}$
be a bounded complex of 
$\Lambda$-modules
on $X$.
Assume that
$f\colon X\to S$
is locally acyclic relatively to ${\cal F}$.

{\rm 1.}
Let ${\cal G}$
be a complex bounded above of 
$\Lambda$-modules
on $S$.
Then,
the canonical morphism
{\rm (\ref{eqPs})} induces an isomorphism
\begin{equation}
p_1^*{\cal F}\otimes^L_\Lambda
p_2^*{\cal G}
\to 
R\Psi_f({\cal F}\otimes^L_\Lambda
f^*{\cal G})
\label{eqla1}
\end{equation}
on $X\overset\gets\times_SS$.

{\rm 2.}
Let $h\colon W\to X$ be
a morphism of schemes
and let ${\cal G}$
be a complex bounded above of 
$\Lambda$-modules
on $W$.
Then,
the canonical morphism
{\rm (\ref{eqcan})} defines an
isomorphism
\begin{equation}
p_1^*h^*{\cal F}
\otimes^L_\Lambda
p_1^*{\cal G}
\to
R\overset\gets f_*(
p_2^*{\cal F}
\otimes^L_\Lambda
p_1^*{\cal G})
\label{eqpfv2}
\end{equation}
on $W\overset\gets\times_SS$.
\end{cor}

\proof{
By the assumption 
of local acyclicity and Proposition \ref{lmapp}
(1)$\Rightarrow$(2),
the formation of
$R\Psi_f{\cal F}$
commutes with finite base change
and 
the canonical morphism
$p_1^*{\cal F}\to R\Psi_f{\cal F}$
is an isomorphism.

1.
By Lemma \ref{lmpf1},
{\rm (\ref{eqPs})} induces an isomorphism
(\ref{eqla1}).

2.
By Lemma \ref{lmpf2},
and by the canonical isomorphism
$\overset\gets h^*
p_1^*\to p_1^*h^*$,
the right hand side of
(\ref{eqcan}) is identified with
that of (\ref{eqpfv2}).
Thus, the assertion follows.
\qed}

\medskip
We briefly recall the definition of
additive convolution from \cite[4.1]{TS}.
Let $k$ be a field and
let $A_1={\mathbf A}^1_{(0)}$
and $A_2={\mathbf A}^2_{(0)}$
denote the henselizations 
of the affine line 
and of the affine plane 
at the origins.
Let
$f\colon X\to A_1$
and
$g\colon Y\to A_1$
be morphisms of finite type.
We regard the fiber product
$(X\times Y)_2
=(X\times Y)\times_{A_1\times A_1}A_2$
as a scheme over $A_1$
by the composition of the second
projection and the morphism 
 $a\colon A_2\to A_1$ induced
by the addition $+\colon
{\mathbf A}^2\to {\mathbf A}^1$.
Morphisms of vanishing toposes
$$\begin{CD}
X\overset \gets\times_{A_1}A_1
@<{{\rm pr}_1}<<
(X\times Y)_2\overset \gets\times_{A_2}A_2
@>{{\rm pr}_2}>>
Y\overset \gets\times_{A_1}A_1\\
@.@V{\overset \gets a}VV@.\\
@.
(X\times Y)_2\overset \gets\times_{A_1}A_1
\end{CD}$$
are defined by projections
and by $a\colon A_2\to A_1$.

Let $\Lambda$ be a finite
field of characteristic invertible in $k$.
For bounded complexes
${\cal F}$ and ${\cal G}$
of $\Lambda$-modules
on $X\overset \gets\times_{A_1}A_1$
and on $Y\overset \gets\times_{A_1}A_1$,
let ${\cal F}\boxtimes{\cal G}$ denote 
${\rm pr}_1^*{\cal F}\otimes
{\rm pr}_2^*{\cal G}$
on
$(X\times Y)_2\overset \gets\times_{A_2}A_2$
and define
the additive convolution
${\cal F}\ast{\cal G}$
on $(X\times Y)_2\overset \gets\times_{A_1}A_1$
by
\begin{equation}
{\cal F}\ast{\cal G}
=
R\overset \gets a_*({\cal F}\boxtimes{\cal G}).
\label{eq*}
\end{equation}

\section{External products}\label{spr}

For the definitions and
basic properties 
of the singular support
of a constructible complex
on a smooth scheme over
a perfect field,
we refer to \cite{Be}
and \cite{CC}.

Let $k$ be a field
and
let $\Lambda$ be a finite
field of characteristic invertible in $k$.

\begin{pr}\label{prpr}
Let $X$ and $Y$ be smooth schemes
over $k$
and ${\cal F}$ and ${\cal G}$
be constructible complexes of
$\Lambda$-modules
on $X$ and on $Y$ respectively.
Assume that
${\cal F}$ 
is micro-supported on
a closed conical subset
$C\subset T^*X$.
Then
${\cal F}\boxtimes^L_\Lambda {\cal G}
={\rm pr}_1^*{\cal F}
\otimes^L_\Lambda
{\rm pr}_2^*{\cal G}$
is micro-supported on
$C\times T^*Y\subset
T^*(X\times Y)$.
\end{pr}

\proof{
It suffices to show that,
for morphisms 
$a\colon W\to X,
b\colon W\to Y,
c\colon W\to Z$ of smooth schemes
over $k$ such that
the pair 
$h=(a,b)\colon W\to X\times Y$
and
$c\colon W\to Z$
is $C\times T^*Y$-transversal,
the morphism
$c\colon W\to Z$
is locally acyclic relatively to $h^*{\cal F}$.
By \cite[Lemma 3.6.9]{CC},
the pair of morphisms
$a\colon W\to X$
and $f=(b,c)\colon W\to Y\times Z$
is $C$-transversal.
Since ${\cal F}$ is assumed 
micro-supported on $C$,
the morphism
$f=(b,c)\colon W\to Y\times Z$
is locally acyclic relatively to $a^*{\cal F}$.

Let $Z'\to Z$ be any finite morphism
and we consider the commutative
diagram
$$\xymatrix{
X&
W'=W\times_ZZ'
\ar[l]_{\!\!\!\!\!\!\!\!\!\!\!\! a'}
\ar[ld]_{h'}
\ar[d]^{b'}
\ar[rd]^{f'}
\ar[r]^{\ \ \ \ \ \ \ \ \ c'}&
Z'\\
X\times Y
\ar[u]\ar[r]&
Y&
Y\times Z'
\ar[u]_r\ar[l]_q}
$$
By Proposition \ref{lmapp}
(3)$\Rightarrow$(1),
it suffices to show that
the canonical morphism
\begin{equation}
p_1^*h^{\prime *}
({\cal F}\boxtimes{\cal G})
\to
R\Psi_{c'}h^{\prime *}
({\cal F}\boxtimes{\cal G})
\label{eqpr1}
\end{equation}
on $W'\overset\gets\times_{Z'}Z'$
is an isomorphism.
For the second term in (\ref{eqpr1}),
we have a canonical isomorphism
\begin{equation}
R\overset\gets r_*
R\Psi_{f'}
(a^{\prime *}{\cal F}
\otimes
f^{\prime *}q^*{\cal G})
\to
R\Psi_{c'}h^{\prime *}
({\cal F}\boxtimes{\cal G}).
\label{eqpr2}
\end{equation}

For the first term in (\ref{eqpr2}),
we apply Corollary \ref{corpf}.1
to $f'\colon W'\to Y\times Z'$
and to $a^{\prime *}{\cal F}$
on $W'$ and $q^*{\cal G}$
on $Y\times Z'$.
Since $f'\colon W'\to Y\times Z'$
is locally acyclic relatively to
$a^{\prime *}{\cal F}$,
the assumption of Corollary \ref{corpf}.1
is satisfied and we obtain a canonical isomorphism
\begin{equation}
p_1^*a^{\prime *}{\cal F}
\otimes
p_2^*q^*{\cal G}
\to
R\Psi_{f'}
(a^{\prime *}{\cal F}
\otimes
f^{\prime *}q^*{\cal G}).
\label{eqpr3}
\end{equation}
Further we apply Corollary \ref{corpf}.2
to $f'\colon W'\to Y\times Z'$
and 
$r\colon Y\times Z'\to Z'$
and to $q^*{\cal G}$
on $Y\times Z'$ and
$a^{\prime *}{\cal F}$
on $W'$.
By the generic local acyclicity
\cite[Corollaire 2.16]{TF},
the second projection
$r\colon Y\times Z'\to Z'$
is locally acyclic relatively to
$q^*{\cal G}$.
Hence the assumption of
Corollary \ref{corpf}.2
is satisfied and we obtain a canonical isomorphism
\begin{equation}
p_1^*a^{\prime *}{\cal F}
\otimes
p_1^*f^{\prime*}q^*{\cal G}
\to
R\overset\gets r_*
(p_1^*a^{\prime *}{\cal F}
\otimes
p_2^*q^*{\cal G})
\label{eqpr4}
\end{equation}
on $W'\overset\gets\times_{Z'}Z'$.
Thus, (\ref{eqpr2})--(\ref{eqpr4})
give an isomorphism
\begin{equation}
p_1^*h^{\prime *}
({\cal F}\boxtimes{\cal G})
=
p_1^*a^{\prime *}{\cal F}
\otimes
p_1^*b^{\prime *}{\cal G}
\to
R\Psi_{c'}h^{\prime *}
({\cal F}\boxtimes{\cal G})
\label{eqpr5}
\end{equation}
and the assertion follows.
\qed}

\medskip

For linear combinations
$A=\sum_am_a\cdot C_a$
and
$A'=\sum_{a'}m'_{a'}\cdot C'_{a'}$
of irreducible components
of closed conical subsets 
$C=\bigcup_aC_a\subset T^*X$
and
$C'=\bigcup_{a'}C'_{a'}\subset T^*Y$
of cotangent bundles,
the external product
$A\boxtimes A'$
is defined by
\begin{equation}
A\boxtimes A'=
\sum_{a,a'}m_am'_{a'}\cdot 
C_a\times C'_{a'}
\label{eqAA'}
\end{equation}
as a linear combination
supported on
$C\times C'\subset
T^*X\times T^*Y=T^*(X\times Y)$.

\begin{thm}\label{thmpr}
Let $X$ and $Y$ be smooth schemes
over a perfect field $k$
and ${\cal F}$ and ${\cal G}$
be constructible complexes of
$\Lambda$-modules.

{\rm 1.}
Assume that
${\cal F}$ and ${\cal G}$
are micro-supported on
closed conical subsets
$C\subset T^*X$ and on
$C'\subset T^*Y$ respectively.
Then
${\cal F}\boxtimes^L_\Lambda {\cal G}$
is micro-supported on
$C\times C'\subset
T^*X\times T^*Y=T^*(X\times Y)$.

{\rm 2.}
We have
\begin{equation}
CC({\cal F}\boxtimes^L_\Lambda {\cal G})
=
CC({\cal F})\boxtimes
CC({\cal G}).
\label{eqCCpr}
\end{equation}

{\rm 3.}
We have
\begin{equation}
SS({\cal F}\boxtimes^L_\Lambda {\cal G})
=
SS({\cal F})\boxtimes
SS({\cal G}).
\label{eqSSpr}
\end{equation}
\end{thm}

We will deduce the assertion 2
from the following multiplicativity
of the Artin conductor
under the convolution.
This is an analogue
for torsion coefficient
of that for 
${\mathbf Q}_\ell$-coefficient
due to Laumon
 \cite[Proposition (2.7.2.1)]{La}.

\begin{lm}[{\rm \cite[Corollary 5.12]{TS}}]\label{lmTS}
Let ${\cal K}$ and ${\cal L}$ be
constructible complexes
of $\Lambda$-modules on 
the strict localization $A_1={\mathbf A}^1_{(0)}$.
Then, for the Artin conductor, we have
\begin{equation}
-a_0({\cal K}*{\cal L})
=
(-a_0{\cal K})\cdot
(-a_0{\cal L}).
\label{eqa*}
\end{equation}
\end{lm}

\proof[Proof of Theorem {\rm\ref{thmpr}}]{
1.
By Proposition \ref{prpr}, 
the external product
${\cal F}\boxtimes^L_\Lambda {\cal G}$
is micro-supported on the intersection
$(C\times T^*Y)\cap 
(T^*X\times C')=C\times C'$.

2.
Write the singular supports
$C=SS({\cal F})
=\bigcup_aC_a$
and 
$C'=SS({\cal G})
=\bigcup_{a'}C'_{a'}$
as the unions of irreducible
components
and set
$CC({\cal F})
=\sum_am_aC_a$
and 
$CC({\cal G})
=\sum_{a'}m_{a'}C'_{a'}$.
Then, by 1,
we have
$CC({\cal F}\boxtimes {\cal G})
=\sum_{a,a'}m_{a,a'}C_a\times C'_{a'}$
for some integers $m_{a,a'}$.
It suffices to show
$m_{b,b'}=m_b\cdot m_{b'}$
for each pair of irreducible components
$C_b$ and $C'_{b'}$.

After shrinking $X$,
we may take a morphism
$f\colon X\to {\mathbf A}^1$
such that
$f$ has an isolated characteristic point $u$,
that $f(u)=0$
and that the section $df$
meets only $C_b$.
Similarly,
after shrinking $Y$,
we may take a morphism
$g\colon Y\to {\mathbf A}^1$
such that
$g$ has an isolated characteristic point $v$,
that $g(v)=0$
and that the section $dg$
meets only $C'_{b'}$.
Let $h\colon X\times Y\to {\mathbf A}^1$
denote the morphism defined by the sum
$f+g$.
Since $dh=df+dg$, the morphism
$h$ has an isolated characteristic point $(u,v)$
with respect to $C\times C'
=\bigcup_{a,a'}C_a\times C'_{a'}$
and that the section $dh$
meets only $C_b\times C'_{b'}$.
Further, we have
$(C_b\times C'_{b'},dh)_{T^*(X\times Y),(u,v)}
=
(C_b,df)_{T^*X,u}
\cdot
(C'_{b'},dg)_{T^*Y,v}\neq0$.
Thus, by the Milnor formula \cite[(5.15)]{CC},
it suffices to show
\begin{equation}
-\dim{\rm tot}
\phi_{(u,v)}({\cal F}\boxtimes
{\cal G},h)
=
(-\dim{\rm tot}
\phi_u({\cal F},f))
\cdot
(-\dim{\rm tot}
\phi_v({\cal G},g)).
\label{eqdt}
\end{equation}

We canonically identify 
$u\overset\gets\times_{{\mathbf A}^1}
{\mathbf A}^1$
with the strict localization
${\mathbf A}^1_{(0)}$.
Then the total dimension
$\dim{\rm tot}
\phi_u({\cal F},f)$
equals the Artin conductor
$a_0\bigl((R\Psi_f{\cal F})|_{u\overset\gets\times_{{\mathbf A}^1}
{\mathbf A}^1}\bigr)$
and similarly for the other terms.
By \cite[Theorem 4.5 (4.5.1)]{TS},
we have an isomorphism
\begin{equation}
R\Psi_f{\cal F}
*
R\Psi_g{\cal G}
\to 
R\Psi_h({\cal F}\boxtimes
{\cal G}).
\label{eqLa}
\end{equation}
The left hand side is the additive convolution
(\ref{eq*}).
Thus, we obtain the equality (\ref{eqdt}) 
by applying Lemma \ref{lmTS} to
${\cal K}=(R\Psi_f{\cal F})|_{u\overset\gets\times_{{\mathbf A}^1}
{\mathbf A}^1}$
and ${\cal L}=(R\Psi_g{\cal G})|_{v\overset\gets\times_{{\mathbf A}^1}
{\mathbf A}^1}$.

3. We may assume ${\cal F}$
and ${\cal G}$ are perverse sheaves.
Then, since the singular support
is the support of the characteristic cycle
by \cite[Proposition 3.19.2]{CC},
the assertion follows from 2.
\qed}

\begin{cor}\label{corpull}
Let $h\colon W\to X$
be a smooth morphism of
smooth schemes over a perfect field $k$.
Then, for a constructible complex ${\cal F}$
of $\Lambda$-modules on $X$,
we have
\begin{equation}
CCh^*{\cal F}=h^!CC{\cal F}.
\label{eqsm}
\end{equation}
\end{cor}

For the definition of the notation $h^!CC{\cal F}$,
we refer to \cite[Definition 5.16]{CC}.

\proof{
Since the assertion is \'etale local on $W$,
we may assume $W=X\times {\mathbf A}^m$
for an integer $m\geqq 0$.
Hence it follows from
Theorem \ref{thmpr}.2.
\qed}
\medskip

In the above proof of Theorem \ref{thmpr}.2,
we deduced (\ref{eqCCpr})
from the multiplicativity
of the Artin conductor
under the convolution
\cite[Corollary 5.12]{TS}.
Conversely, \cite[Corollary 5.12]{TS}
is an immediate consequence of
(\ref{eqCCpr}) 
where
$\dim X=\dim Y=1$.
This crucial case 
can be deduced from
the index formula proved
earlier in \cite[Theorem 3.19]{surface}
and is
essentially equivalent to
\cite[Exemples 2.3.8 (a)]{chi}.
Corollary \ref{corpull} is 
a first step of the proof of
the index formula
\cite[Theorem 7.13]{CC}
in general dimension.
These logical implications
are summarized in the diagram
$$
\begin{matrix}
&&
\left(\text{
\parbox{42mm}
{index formula in dim.\! 2\\
\cite[Theorem 3.19]{surface} or\\
\cite[Exemples 2.3.8 (a)]{chi}}}
\right)\\
&&\Downarrow\\
\left(\text{
\parbox{38mm}
{multiplicativity of\\
the Artin conductor\\
\cite[Corollary 5.12]{TS}}}
\right)
&\Leftrightarrow&
\left(\text{
\parbox{38mm}
{multiplicativity of\\
characteristic cycles\\
(\ref{eqCCpr}) in dim.\! $1+1$}}
\right)\\
\Downarrow&&\\
\left(\text{
\parbox{38mm}
{compatibility with\\
smooth pull-back\\
Corollary \ref{corpull}}}
\right)
&\Rightarrow&
\left(\text{
\parbox{45mm}
{index formula in dim$\geqq 2$\\
\cite[Theorem 7.13]{CC}}}
\right).
\end{matrix}
$$

\begin{cor}\label{corpr}
Let ${\cal F}$ 
and ${\cal G}$ be constructible complexes
of $\Lambda$-modules 
on  a smooth scheme $X$ over $k$.
Assume that the intersection
$SS({\cal F})\cap SS({\cal G})
\subset T^*X$
of the singular supports
is a subset of the $0$-section
$T^*_XX\subset T^*X$.
Then, the canonical morphism
\begin{equation}
{\cal G}\otimes^L
R{\cal H}om_X({\cal F},\Lambda)
\to
R{\cal H}om_X({\cal F},{\cal G})
\label{eqRHom}
\end{equation}
is an isomorphism.
\end{cor}

\proof{
Set
$C=SS({\cal F})$ and $C'=SS({\cal G})$.
The assumption
$C\cap C'\subset T^*_XX$
implies that the
diagonal $\delta\colon X\to X\times X$
is $C'\times C$-transversal.
Since $SSD_X{\cal F}=SS{\cal F}$
by \cite[Corollary 2.27]{CC},
the external product
${\cal G}\boxtimes D_X{\cal F}$
is micro-supported on $C'\times C
\subset T^*X\times T^*X
=T^*(X\times X)$
by Theorem \ref{thmpr}.1.
Since the canonical morphism
${\cal G}\boxtimes D_X{\cal F}
\to
R{\cal H}om_{X\times X}({\rm pr}_2^*{\cal F},
{\rm pr}_1^!{\cal G})$
is an isomorphism
by \cite[(3.1.1)]{SGA5},
the diagonal $\delta\colon X\to X\times X$
is $R{\cal H}om_{X\times X}({\rm pr}_2^*{\cal F},
{\rm pr}_1^!{\cal G})$-transversal
by \cite[Proposition 5.6]{CC}.
Thus, the assertion follows from
\cite[Proposition 5.3.2 (1)$\Rightarrow$(2)]{CC}.
\qed}
\medskip

Recall that a closed subset
of a vector bundle
said to be {\em conical} if it is
stable under the action of
the multiplicative group ${\mathbf G}_m$.

\begin{df}\label{dftrans}
Let $f\colon X\to Y$ be 
a morphism
of smooth schemes over $k$
Let $C\subset T^*X$ 
and $C'\subset T^*Y$ be
closed conical subsets
of the cotangent bundles.

For $x\in X$,
we say that $f\colon X\to Y$
is $(C,C')$-{\em transversal}
if for $\omega \in T^*Y\times_Yy$
at $y=f(x)\in Y$,
the conditions
$\omega\in C',f^*\omega\in C$
imply $\omega=0$.
We say 
$f\colon X\to Y$
is $(C,C')$-{\em transversal}
if $f\colon X\to Y$
is $(C,C')$-transversal at every $x\in X$.
Or equivalently,
if the intersection
$df^{-1}(C)\cap f^*C'$
of the inverse image of
$C$ by $df\colon 
X\times_YT^*Y\to T^*X$
and $f^*C'=X\times_YC'$
in $X\times_YT^*Y$
is a subset of the $0$-section.
\end{df}

\begin{lm}\label{lmCC'}
Let $f\colon X\to Y$
be a morphism of
smooth schemes over $k$
and 
let $\gamma\colon X\to X\times Y$
be the graph of $f$.
For closed conical subsets
$C\subset T^*X$ and
$C'\subset T^*Y$
the following conditions are
equivalent:

{\rm (1)}
$f$ is $(C,C')$-transversal.

{\rm (2)}
$\gamma$ is 
$C\times C'$-transversal.

Further, if the condition {\rm (2)}
is satisfied,
the closed subset
$\gamma^o(C\times C')\subset T^*X$
equals the subset $C+f^*C'\subset T^*X$
consisting of the sum $\alpha+\beta$
of $\alpha\in C$ and $\beta \in f^*C'$.
\end{lm}

\proof{
The condition (1) 
is equivalent to the following the condition:

(1$'$)
For $\beta\in f^*C'$,
if $df(\beta)\in C$,
then we have $\beta=0$.

\noindent
The condition (2) 
is equivalent to the following the condition:

(2$'$)
For $\alpha\in C$ and $\beta\in f^*C'$,
if $\alpha+df(\beta)=0$,
then we have $\alpha=0$ and $\beta=0$.

\noindent
Hence the conditions (1) and (2)
are equivalent.
Since $\gamma^*(C\times C')=
C\times_Xf^*C'$,
we obtain
$\gamma^o=C+f^*C'$.
\qed}

\medskip
For a separated morphism
$h\colon W\to X$ of finite
type and for a constructible
complex ${\cal F}$ of 
$\Lambda$-modules,
a canonical morphism
$c_{h,{\cal F}}
\colon 
h^*{\cal F}\otimes Rh^!\Lambda
\to Rh^!{\cal F}$
is defined in \cite[(8.13)]{CC}
as the adjoint of
the morphism
$Rh_!(h^*{\cal F}\otimes Rh^!\Lambda)
\overset \sim \gets
{\cal F}\otimes Rh_!Rh^!\Lambda
\to {\cal F}$
induced by the adjunction
$Rh_!Rh^!\Lambda\to \Lambda$.

\begin{pr}\label{prFG}
Let $f\colon X\to Y$
be a morphism of
smooth schemes over $k$
and 
let $\gamma\colon X\to X\times Y$
be the graph of $f$.
Let ${\cal F}$ be a constructible complex
of $\Lambda$-modules 
on $X$ and set $C=SS({\cal F})$.

Let ${\cal G}$ be a constructible complex
of $\Lambda$-modules 
on $Y$ and set $C'=SS({\cal G})$.
Assume $f$ is $(C,C')$-transversal.
Then, the canonical morphism
\begin{equation}
\gamma^*({\cal F}\boxtimes {\cal G})
\otimes R\gamma^!\Lambda
\to
R\gamma^!({\cal F}\boxtimes {\cal G})
\label{eqgm}
\end{equation}
is an isomorphism
and
${\cal F}\otimes f^*{\cal G}=
\gamma^*({\cal F}\boxtimes {\cal G})$
is micro-supported
on $C+f^*C'\subset T^*X$
consisting of the sum $\alpha+\beta$
of $\alpha\in C$ and $\beta \in f^*C'$.
\end{pr}

The morphism (\ref{eqgm})
is the same as \cite[(5.3)]{BG}.

\proof{
The assumption that
$f$ is $(C,C')$-transversal
means that
$\gamma$ is 
$C\times C'$-transversal
by Lemma \ref{lmCC'}.
Since ${\cal F}\boxtimes {\cal G}$
is micro-supported on $C\times C'$
by Theorem \ref{thmpr}.1,
the morphism $\gamma$ is 
${\cal F}\boxtimes {\cal G}$-transversal
by \cite[Proposition 5.6 (1)$\Rightarrow$(2)]{CC}.
Thus, the morphism
(\ref{eqgm}) is an isomorphism.
Further,
${\cal F}\otimes f^*{\cal G}=
\gamma^*({\cal F}\boxtimes {\cal G})$
is micro-supported
on $\gamma^o(C\times C')=C+f^*C'$
by \cite[Lemma 2.11.4 (1)$\Rightarrow$(2)]{CC} and Lemma \ref{lmCC'}.
\qed}

\begin{cor}\label{corFG}
Let the notation be as in
Proposition {\rm \ref{prFG}} and
assume that $f$ is $C$-transversal.
Then, 
for every constructible complex ${\cal G}$
of $\Lambda$-modules 
on $Y$,
the conclusion of Proposition {\rm \ref{prFG}}
is satisfied.
\end{cor}

The conclusion of Corollary \ref{corFG}
is shown to be equivalent to the 
local acyclicity of $f$ in
\cite[Theorem B.2]{BG}.

\proof{
Since $f$ is $(C,T^*Y)$-transversal,
it is $(C,C')$-transversal
for any closed conical subset
$C'\subset T^*Y$.
\qed}

\end{document}